\documentclass[12pt, preprint, times]{elsarticle}

\usepackage{amsmath}
\usepackage{array}
\usepackage{amssymb}
\usepackage{amsthm}
\usepackage{mathtools}
\mathtoolsset{showonlyrefs}

\newtheorem{theorem}{Theorem}

\newtheorem{lemma}{Lemma}

\theoremstyle{definition}
\newtheorem*{remark}{Remark}
\newtheorem{definition}{Definition}

\numberwithin{equation}{section}

\begin{document}

\begin{frontmatter}

\title{\fontsize{16.5}{22}\selectfont A Threshold-Crossing Condition for the Least Colossally Abundant Counterexample to Robin's Inequality}

\author[CRI]{Bruce Zimov}
\address[CRI]{Calimesa Research Institute, Calimesa, CA, USA} 
\ead{bzimov@calri.org}

\begin{abstract}
We study the constraints that a hypothetical least Colossally Abundant (CA) counterexample to Robin's Inequality would place on the non-trivial zeros of $\zeta(s)$. In the single-prime case we express the least CA counterexample $N$ as $M\cdot p$, where $M$ is its CA predecessor and $p$ is prime. Using the von Mangoldt explicit formula for $\psi(x)$, we express $\ln p$ in terms of a sum over the non-trivial zeros of $\zeta(s)$ and split that expression into its baseline contribution $\mathcal{R}_0(p)$ and its off-line distortion. When both converge, we show that the off-line distortion is negative whenever $1+1\big/\big(e^{\mathcal{R}_0(p)}+e^{2\mathcal{R}_0(p)}+\cdots+e^{(v+1)\mathcal{R}_0(p)}\big)\le\ln(\ln M+\mathcal{R}_0(p))\big/\ln\ln M$, where $v$ is the $p$-adic valuation of $M$. In the case where the least CA counterexample is the product of its CA predecessor and two primes (the semiprime case), the logarithm of the product of the two primes is the sum of two von Mangoldt expansions. An analogous inequality implies that at least one of the off-line distortions of the two expansions is negative.
\end{abstract}
\begin{keyword}
Colossally Abundant Numbers \sep Robin's Inequality \sep Riemann~Hypothesis
\\
MSC codes: 11A25 \sep 11M26 \sep 11N37 \sep 11N56
\end{keyword}

\end{frontmatter}

\section{Introduction}
\label{sec:intro}
In 1984, Guy Robin \cite{RefF} established a remarkable equivalence, demonstrating that the Riemann Hypothesis is true if and only if a specific inequality, now known as Robin's Inequality, holds for all integers $n>5040$. This inequality, given by $G(n)<e^\gamma$, where $\gamma$ is the Euler--Mascheroni constant, has since become a central focus for researchers. 
\begin{samepage}
\begin{definition}
\label{def:G}
For an integer $n\ge2$,
\begin{equation}
G(n) \coloneqq \displaystyle \frac{\sigma(n)}{n \ln\ln n},\end{equation} where $\sigma(n)$ is the sum of divisors function.
\end{definition}
\end{samepage}
A key result in this area, proven by Akbary and Friggstad \cite{RefA}, shows that if a counterexample to the inequality exists, then the least such counterexample must be a superabundant number. These numbers are positive integers $n$ defined by the property that the ratio $\displaystyle\frac{\sigma(n)}{n}$ is greater than that for any preceding integer. A critical subset of superabundant numbers, known as colossally abundant (CA) numbers, provides a more structured and manageable set. These were first studied by Alaoglu and Erdős in 1944 \cite{RefB} and by Erdős and Nicolas in 1975 \cite{RefD}.  The existence of a least CA counterexample is a question that, if answered in the negative, would provide a powerful, unconditional proof of the Riemann Hypothesis \cite[Lemma 7.1]{RefC}. Resolving the Riemann Hypothesis itself is not the goal of this paper. Our goal is to determine the impact that a least CA counterexample would have on the non-trivial zeros of $\zeta(s)$.

Let $N$ be the least CA counterexample to Robin's Inequality such that $G(N) \ge e^\gamma$, and let $M$ be its immediate CA predecessor satisfying $G(M) < e^\gamma$. By Morrill and Platt \cite[Theorem 5]{RefP}, $N$ is far greater than 55440, the CA number following 5040. Hence $M>5040$. By Broughan's Lemma~6.15 \cite{RefC}, $N$ and $M$ are related by $N = M\cdot p$ with $p$ prime or $N = M \cdot pq$ with $p<q$ distinct primes. In the single-prime case, $N = M\cdot p$ with $v := \nu_p(M) \ge 0$ the $p$-adic valuation of $M$. Either $p$ is a \emph{new prime} not already dividing $M$ in which case $v=0$ and $p^{1}$ is the prime-power content of $N$, or $p$ is an \emph{interior prime} already dividing $M$ in which case $v \ge 1$ and $p^{v+1}$ is the prime-power content of $N$. For every CA number $M>5040$, the exponent of each prime stays the same or increases along the CA sequence \cite[proof of Lemma 6.15]{RefC}. Since $\nu_2(5040)=4$, this gives $\nu_2(M)\ge4$. In particular, $p^{v+1}>2$ for every prime $p$ and $v=\nu_p(M)$.

Let $\rho=\beta_\rho+i\gamma_\rho=\frac{1}{2}+\delta_\rho+i\gamma_\rho$ denote a non-trivial zero of $\zeta(s)$, where $0<\beta_\rho<1$ and $\delta_\rho$ measures its horizontal displacement from the critical line.

\begin{definition}[Critical-line and off-line addends]
\label{def:addends}
For a sum over the non-trivial zeros $\rho$ of $\zeta(s)$, its \emph{critical-line addend} is the sum of its terms over the zeros with $\mathrm{Re}(\rho)=\frac{1}{2}$, and its \emph{off-line addend}, also called its \emph{off-line distortion}, is the sum of its terms over the zeros with $\mathrm{Re}(\rho)\neq\frac{1}{2}$.
\end{definition}

\noindent Sums over the non-trivial zeros are taken in the symmetric sense $\lim_{T\to\infty}\sum_{|\gamma_\rho|<T}$, as in Davenport \cite[Chapter 17]{RefE}, whose argument shows that the full sum converges. Pairing each zero with its conjugate turns such a sum into a sum over $\gamma_\rho>0$ of the real terms $x^\rho/\rho+x^{\bar\rho}/\bar\rho$. Because conjugate zeros share the same real part, no pair is split between the critical-line and off-line addends.

By expressing $\ln p$ directly via von Mangoldt's explicit formula, we find that the sum over the non-trivial zeros of $\zeta(s)$ allows us to isolate and analyze how off-line zero distortions behave in the case of the least CA counterexample. Our main results treat the case in which both the critical-line addend and the off-line addend of this sum converge.

Alaoglu and Erdős \cite{RefB} conjectured that $p^x$ and $q^x$ are both rational, for primes $p\ne q$, only if $x$ is an integer. This conjecture would imply that the quotient of consecutive CA numbers is prime, and it follows from the four exponentials conjecture. Because the four exponentials conjecture has not been proven (see Lagarias \cite{RefL}, Lang \cite{RefN}), we extend the analysis to the semiprime case.

\section{Background Results}
\begin{definition}[Broughan, 2017]\label{def:CA_properties}
A positive integer $n$ is a CA number if there exists an exponent $\epsilon > 0$ such that $\displaystyle\frac{\sigma(k)}{k^{1+\epsilon}}$ reaches its maximum at $n$, where $\sigma(k)$ is the sum of divisors function. For all positive integers $k$,
\[
\displaystyle\frac{\sigma(k)}{k^{1+\epsilon}}\le \displaystyle\frac{\sigma(n)}{n^{1+\epsilon}}.
\]
\end{definition}
\begin{theorem}
\label{7.1}
If a counterexample to Robin's Inequality exists for some integer $n > 5040$, then a counterexample which is a CA number also exists.
\begin{proof}
This is Broughan \cite[Lemma 7.1]{RefC}.
\end{proof}
\end{theorem}
\begin{theorem}
\label{6.15}
The quotient of two consecutive CA numbers is either a prime or the product of two distinct primes.
\begin{proof}
This is Broughan \cite[Lemma 6.15]{RefC}.
\end{proof}
\end{theorem}
\begin{theorem}
\label{super}
The smallest integer $n>5040$ which does not satisfy Robin's Inequality must be a superabundant number.
\begin{proof}
This is Akbary and Friggstad \cite[Theorem 3]{RefA}.
\end{proof}
\end{theorem}
\begin{samepage}
\begin{theorem}
\label{thm:MorrillPlatt}
Robin's Inequality holds for all integers $n$ with $5040<n\le10^{10^{13.11485}}$.
\begin{proof}
This is Morrill and Platt \cite[Theorem 5]{RefP}.
\end{proof}
\end{theorem}
\end{samepage}
\begin{remark}
The least CA counterexample to Robin's Inequality is not necessarily the least counterexample to Robin's Inequality, which could be superabundant but not CA. However, by Theorem~\ref{7.1}, if one could prove that no CA number greater than 5040 is a counterexample to Robin's Inequality, then Robin's Inequality must hold for all integers $n > 5040$.
\end{remark}
\begin{theorem}[von Mangoldt Explicit Formula]
\label{thm:von_mangoldt_explicit}
Let $\Lambda(n)$ be the von Mangoldt function, $\Lambda(n)=\ln q$ if $n=q^k$ for a prime $q$ and an integer $k\ge1$, and $\Lambda(n)=0$ otherwise. Let $\psi(x)=\sum_{n\le x}\Lambda(n)$ be Chebyshev's function. Since $\psi(x)$ jumps at each prime power, let $\psi_0(x)=\psi(x)$ if $x$ is not a prime power, and $\psi_0(x)=\psi(x)-\tfrac12\Lambda(x)$ if it is. Then for any real $x > 1$
\begin{equation}
\psi_0(x) = x - \sum_{\rho} \frac{x^\rho}{\rho} - \ln(2\pi) - \frac{1}{2}\ln(1 - x^{-2}),
\end{equation}
where the sum over the non-trivial zeros $\rho$ of $\zeta(s)$ is to be understood in the symmetric sense as
\begin{equation}
\lim_{T \to \infty} \sum_{|\gamma_\rho| < T} \frac{x^\rho}{\rho}.
\end{equation}
Moreover, at every positive integer $n\ge2$,
\begin{equation}
\psi(n)-\psi(n-1)=\Lambda(n),
\end{equation}
and at every integer $x\ge2$,
\begin{equation}
\psi(x) = \psi_0(x) + \frac{1}{2}\Lambda(x).
\end{equation}
\end{theorem}
\begin{proof}
The definitions of $\Lambda(n)$ and $\psi(x)$ are those of Davenport \cite[Chapters 7 and 17]{RefE}. The explicit formula, with $\psi_0(x)$ as defined here, is in Davenport \cite[Chapter 17]{RefE}. Broughan \cite[Theorem 2.1]{RefC} also states it, for $x$ not a prime power, and records that at a prime power the right-hand side is the average of $\psi(x-)$ and $\psi(x+)$. The last two identities follow from the definitions. First, by the definition of $\psi$, $\psi(n)=\Lambda(1)+\cdots+\Lambda(n)$ and $\psi(n-1)=\Lambda(1)+\cdots+\Lambda(n-1)$, so $\psi(n)-\psi(n-1)=\Lambda(n)$. Second, if $x$ is a prime power, the last identity is the definition of $\psi_0(x)$, and otherwise $\Lambda(x)=0$ and $\psi_0(x)=\psi(x)$.
\end{proof}

\section{Main Results}

\subsection{Spectral Identity for $\ln p$}

\begin{lemma}[Spectral Identity for $\ln p$]\label{lem:ln-p-identity}
Let $p$ be a prime and $v\ge0$ an integer with $p^{v+1}>2$. Then
\begin{equation}
\begin{alignedat}{2}
\ln p &= \big(2 && -\Lambda(p^{v+1}-1)\big) \\
&&& - \sum_{\gamma_\rho > 0} \frac{4}{|\rho|} \Bigg[ (p^{v+1})^{\frac{1}{2}+\delta_\rho} \cos\left(\gamma_\rho \ln p^{v+1} - \phi(\delta_\rho,\gamma_\rho)\right) \\
&&& - (p^{v+1}-1)^{\frac{1}{2}+\delta_\rho} \cos\left(\gamma_\rho \ln(p^{v+1}-1) - \phi(\delta_\rho,\gamma_\rho)\right) \Bigg] \\
&&& - \ln \left( \frac{1 - (p^{v+1})^{-2}}{1 - (p^{v+1}-1)^{-2}} \right),
\end{alignedat}
\label{eq:ln_p_phase_form}
\end{equation}
where $\Lambda(n)$ is the von Mangoldt function of Theorem~\ref{thm:von_mangoldt_explicit} and the phase angle for $\gamma_\rho > 0$ is defined by
\begin{equation}
\phi(\delta_\rho,\gamma_\rho) \coloneqq \arctan\left( \frac{\gamma_\rho}{\frac{1}{2} + \delta_\rho} \right).
\end{equation}

\end{lemma}
\begin{proof}

Since $p^{v+1}$ is a power of the prime $p$, $\Lambda(p^{v+1})=\ln p$ by definition, and by Theorem~\ref{thm:von_mangoldt_explicit} at $n=p^{v+1}$ the same value is the difference of $\psi(p^{v+1})$ and $\psi(p^{v+1}-1)$. Hence
\begin{equation}
\ln p = \Lambda\left(p^{v+1}\right) = \psi\left(p^{v+1}\right) - \psi\left(p^{v+1}-1\right).
\label{eq:final_section1_identity}
\end{equation}

To relate this step difference to the non-trivial zeros $\rho = \beta_\rho + i\gamma_\rho$ of $\zeta(s)$, we apply the explicit formula and the last identity of Theorem~\ref{thm:von_mangoldt_explicit} at each of the two boundaries.

Both boundaries satisfy the $x>1$ requirement of Theorem~\ref{thm:von_mangoldt_explicit}, since $p^{v+1}>2$. The relation at $x = p^{v+1}$ and $x = p^{v+1}-1$ reads
\begin{align}
\psi(p^{v+1}) &= \psi_0(p^{v+1}) + \frac{1}{2}\ln p, \label{eq:psi_upper} \\
\psi(p^{v+1}-1) &= \psi_0(p^{v+1}-1) + \frac{1}{2}\Lambda(p^{v+1}-1). \label{eq:psi_lower}
\end{align}

Substituting \eqref{eq:psi_upper} and \eqref{eq:psi_lower} into the step-difference equation \eqref{eq:final_section1_identity} gives
\begin{equation}
\ln p = \left( \psi_0(p^{v+1}) + \frac{1}{2}\ln p \right) - \left(\psi_0(p^{v+1}-1) + \frac{1}{2}\Lambda(p^{v+1}-1)\right).
\end{equation}

Subtracting $\frac{1}{2}\ln p$ from both sides and adding $\frac{1}{2}\Lambda(p^{v+1}-1)$ yields
\begin{equation}
\frac{1}{2}\ln p + \frac{1}{2}\Lambda(p^{v+1}-1) = \psi_0(p^{v+1}) - \psi_0(p^{v+1}-1).
\label{eq:psi_0_diff}
\end{equation}

Theorem~\ref{thm:von_mangoldt_explicit} at $p^{v+1}$ and at $p^{v+1}-1$ gives
\begin{equation}
\psi_0(p^{v+1}) = p^{v+1} - \sum_{\rho} \frac{(p^{v+1})^\rho}{\rho} - \ln(2\pi) - \frac{1}{2}\ln\big(1 - (p^{v+1})^{-2}\big),
\end{equation}
\begin{equation}
\psi_0(p^{v+1}-1) = (p^{v+1}-1) - \sum_{\rho} \frac{(p^{v+1}-1)^\rho}{\rho} - \ln(2\pi) - \frac{1}{2}\ln\big(1 - (p^{v+1}-1)^{-2}\big).
\end{equation}
Each sum over the non-trivial zeros converges, so the difference of the two sums is the sum of the termwise differences. Subtracting the second display from the first gives
\begin{equation}
\psi_0(p^{v+1}) - \psi_0(p^{v+1}-1) = 1 - \sum_{\rho} \frac{(p^{v+1})^\rho - (p^{v+1}-1)^\rho}{\rho} - \frac{1}{2}\ln\left( \frac{1 - (p^{v+1})^{-2}}{1 - (p^{v+1}-1)^{-2}} \right).
\end{equation}

Substituting this expansion into \eqref{eq:psi_0_diff} and multiplying through by $2$ establishes the spectral identity for $\ln p$
\begin{equation}
\ln p = \big(2-\Lambda(p^{v+1}-1)\big) - 2 \sum_{\rho} \frac{(p^{v+1})^\rho - (p^{v+1}-1)^\rho}{\rho} - \ln\left( \frac{1 - (p^{v+1})^{-2}}{1 - (p^{v+1}-1)^{-2}} \right).
\label{eq:exact_ln_p_limit}
\end{equation}

Because non-trivial zeros occur in complex conjugate pairs $\rho = \beta_\rho + i\gamma_\rho$ and $\bar{\rho} = \beta_\rho - i\gamma_\rho$ with $\gamma_\rho > 0$, we group the sum over upper half-plane zeros $\gamma_\rho > 0$. For any real $x > 1$,
\begin{equation}
\begin{split}
\frac{x^\rho}{\rho} + \frac{x^{\bar{\rho}}}{\bar{\rho}} &= \frac{x^{\beta_\rho + i\gamma_\rho}(\beta_\rho - i\gamma_\rho) + x^{\beta_\rho - i\gamma_\rho}(\beta_\rho + i\gamma_\rho)}{\beta_\rho^2 + \gamma_\rho^2} \\
&= \frac{2 x^{\frac{1}{2}+\delta_\rho}}{\left(\frac{1}{2}+\delta_\rho\right)^2 + \gamma_\rho^2} \left[ \left(\frac{1}{2}+\delta_\rho\right)\cos(\gamma_\rho \ln x) + \gamma_\rho \sin(\gamma_\rho \ln x) \right].
\end{split}
\label{eq:pair_expansion}
\end{equation}
Applying \eqref{eq:pair_expansion} to the difference terms $x = p^{v+1}$ and $x = p^{v+1}-1$ converts \eqref{eq:exact_ln_p_limit} into a real-valued trigonometric sum
\begin{equation}
\begin{split}
\ln p = \big(2-\Lambda(p^{v+1}-1)\big) - \sum_{\gamma_\rho > 0} \frac{4}{\left(\frac{1}{2}+\delta_\rho\right)^2 + \gamma_\rho^2} \Bigg[ &\left(\frac{1}{2}+\delta_\rho\right) \Big( (p^{v+1})^{\frac{1}{2}+\delta_\rho} \cos\left(\gamma_\rho \ln p^{v+1}\right) \\
&- (p^{v+1}-1)^{\frac{1}{2}+\delta_\rho} \cos\left(\gamma_\rho \ln(p^{v+1}-1)\right) \Big) \\
&+\gamma_\rho \Big( (p^{v+1})^{\frac{1}{2}+\delta_\rho} \sin\left(\gamma_\rho \ln p^{v+1}\right) \\
&- (p^{v+1}-1)^{\frac{1}{2}+\delta_\rho} \sin\left(\gamma_\rho \ln(p^{v+1}-1)\right) \Big) \Bigg] \\
&- \ln \left( \frac{1 - (p^{v+1})^{-2}}{1 - (p^{v+1}-1)^{-2}} \right).
\end{split}
\end{equation}

Applying the harmonic addition identity, valid for $u>0$ and $w>0$,
\[
u \cos\theta + w \sin\theta = \sqrt{u^2 + w^2}\, \cos\left(\theta - \arctan(w/u)\right)
\]
with $u = \frac{1}{2}+\delta_\rho = \beta_\rho > 0$ and $w = \gamma_\rho > 0$, the phase becomes $\phi(\delta_\rho,\gamma_\rho)$ and the modulus simplifies to
\[
|\rho| = \sqrt{\left(\frac{1}{2}+\delta_\rho\right)^2 + \gamma_\rho^2}.
\]
\emergencystretch=3em Canceling one power of $|\rho|$ in the denominator yields the identity \eqref{eq:ln_p_phase_form} of Lemma~\ref{lem:ln-p-identity}.
\end{proof}

\subsection{Zero Spectrum Partition and Baseline Definition}
\label{sec:baseline}
Lemma~\ref{lem:ln-p-identity} concerns the full sum over the non-trivial zeros. Convergence of the full sum does not give convergence of either addend, because a conditionally convergent sum need not converge after being split. Splitting the sum in the identity of Lemma~\ref{lem:ln-p-identity} into its critical-line and off-line addends, we define the baseline contribution $\mathcal{R}_0(p)$ as everything in that identity apart from the off-line addend, namely the leading constant $2-\Lambda(p^{v+1}-1)$, the log-correction term, and the critical-line addend together

\begin{equation}
\begin{alignedat}{2}
\mathcal{R}_0(p) &\coloneqq \big(2 && -\Lambda(p^{v+1}-1)\big) \\
&&& - \sum_{\substack{\gamma_\rho > 0 \\ \mathrm{Re}(\rho) = \frac{1}{2}}} \frac{4}{|\rho_0|} \Bigg[ \sqrt{p^{v+1}} \cos\left(\gamma_\rho \ln p^{v+1} - \phi(0,\gamma_\rho)\right) \\
&&& - \sqrt{p^{v+1}-1} \cos\left(\gamma_\rho \ln(p^{v+1}-1) - \phi(0,\gamma_\rho)\right) \Bigg] \\
&&& - \ln \left( \frac{1 - (p^{v+1})^{-2}}{1 - (p^{v+1}-1)^{-2}} \right),
\end{alignedat}
\label{eq:R0-def}
\end{equation}
\noindent where $|\rho_0| = \sqrt{\frac{1}{4} + \gamma_\rho^2}$ and $\phi(0,\gamma_\rho) = \arctan(2\gamma_\rho)$. By the definition of $\Lambda$ in Theorem~\ref{thm:von_mangoldt_explicit}, $\Lambda(p^{v+1}-1)=0$ unless $p^{v+1}-1$ is a prime power.

Each addend of the sum in Lemma~\ref{lem:ln-p-identity} has a partial sum over $0<\gamma_\rho<T$ at every $T$. These partial sums depend on the pair $p$, $v$. This paper considers them only for the primes dividing $N/M$, where $N$ is the least CA counterexample, $M$ its CA predecessor, each prime with $v$ the $p$-adic valuation of $M$. In the \emph{convergent case} the partial sums of both addends converge to real numbers as $T\to\infty$ for each such prime. In the \emph{divergent case} the partial sums of at least one addend fail to converge to real numbers as $T\to\infty$ for at least one such prime, whether by tending to $-\infty$, by tending to $+\infty$, or by oscillating. The main results of this paper treat only the convergent case.

\subsection{Isolation of the Off-Line Distortion Term $\Delta(p, \{\delta_\rho\})$}
\label{sec:offline-distortion}
Correspondingly, we define the off-line distortion $\Delta(p, \{\delta_\rho\})$ as the off-line addend of the same identity
\begin{equation}
\begin{split}
\Delta(p, \{\delta_\rho\}) \coloneqq -\sum_{\substack{\gamma_\rho > 0 \\ \mathrm{Re}(\rho) \neq \frac{1}{2}}} \frac{4}{|\rho|} \Bigg[ &(p^{v+1})^{\frac{1}{2}+\delta_\rho} \cos\left(\gamma_\rho \ln p^{v+1} - \phi(\delta_\rho,\gamma_\rho)\right) \\
&- (p^{v+1}-1)^{\frac{1}{2}+\delta_\rho} \cos\left(\gamma_\rho \ln(p^{v+1}-1) - \phi(\delta_\rho,\gamma_\rho)\right) \Bigg].
\end{split}
\end{equation}
By construction,
\begin{equation}
\label{eq:R0-Delta-sum}
\mathcal{R}_0(p) + \Delta(p, \{\delta_\rho\}) = \ln p, \qquad\text{equivalently}\qquad \Delta(p, \{\delta_\rho\}) = \ln p - \mathcal{R}_0(p).
\end{equation}

\begin{lemma}[Master Tracking Inequality]\label{lem:master-tracking}
Let $N=M\cdot p$ be the least CA counterexample to Robin's Inequality
in the single-prime case, with $M$ its CA predecessor and $v$ the $p$-adic valuation of $M$. Then
\begin{equation}
1+\frac{1}{p+p^2+\cdots+p^{v+1}} > \frac{\ln(\ln M+\ln p)}{\ln\ln M}.
\end{equation}
\end{lemma}

\begin{proof}

Since $N$ is the least CA counterexample, $G(N)\ge e^\gamma>G(M)$, hence $G(N)>G(M)$.

The abundancy index scales as
\begin{equation}
\frac{\sigma(N)}{N} = \frac{\sigma(M)}{M} \cdot \left[ \frac{p^{v+2} - 1}{p(p^{v+1} - 1)} \right].
\label{eq:abundancy-scaling}
\end{equation}

For $p>1$,
\begin{equation}
\frac{p^{v+2}-1}{p(p^{v+1}-1)}=\frac{p^{v+2}-1}{p^{v+2}-p}=\frac{1-p^{-(v+2)}}{1-p^{-(v+1)}}=1+\frac{1}{p+p^2+\cdots+p^{v+1}}.
\label{eq:abundancy-factor}
\end{equation}
The first two equalities expand the denominator and divide its numerator and denominator by $p^{v+2}$. The last holds because $p+p^2+\cdots+p^{v+1}=p(p^{v+1}-1)/(p-1)$. Hence
\begin{equation}
\frac{\sigma(N)}{N} = \frac{\sigma(M)}{M}\cdot\left(1+\frac{1}{p+p^2+\cdots+p^{v+1}}\right).
\label{eq:sigma-identity}
\end{equation}

Since $\ln N=\ln M+\ln p$, we have $\ln\ln N=\ln(\ln M+\ln p)$. Therefore
\begin{equation}
\begin{split}
G(N) &= \frac{\sigma(N)/N}{\ln\ln N} \\
&= \frac{\sigma(M)/M}{\ln\ln M}\cdot\frac{1+\dfrac{1}{p+p^2+\cdots+p^{v+1}}}{\ln(\ln M+\ln p)\big/\ln\ln M} \\
&= G(M)\cdot\frac{1+\dfrac{1}{p+p^2+\cdots+p^{v+1}}}{\ln(\ln M+\ln p)\big/\ln\ln M}.
\end{split}
\end{equation}

Since $M>5040$, we have $\ln\ln M>0$, and since also $p\ge2$, $\ln M+\ln p>1$.
Hence $\ln(\ln M+\ln p)>0$, so the denominator of the last fraction is positive. As $G(N)>G(M)>0$, dividing by $G(M)$ shows that the last fraction exceeds $1$, which is the inequality of the lemma.
\end{proof}

Let $\Delta=\Delta(p,\{\delta_\rho\})$ denote the off-line distortion, such that $\ln p=\mathcal{R}_0(p)+\Delta$ by \eqref{eq:R0-Delta-sum}. Then $p=e^{\mathcal{R}_0(p)+\Delta}$. Define
\begin{equation}
\mathbf{A}(\mathcal{R}_0(p),\Delta) \coloneqq 1+\frac{1}{\sum_{k=1}^{v+1}e^{k(\mathcal{R}_0(p)+\Delta)}},
\qquad
\mathbf{B}(\mathcal{R}_0(p),\Delta) \coloneqq \frac{\ln(\ln M+\mathcal{R}_0(p)+\Delta)}{\ln\ln M}.
\label{eq:AB-definition}
\end{equation}
Since $p=e^{\mathcal{R}_0(p)+\Delta}$, the inequality of Lemma~\ref{lem:master-tracking} reads
\begin{equation}
\mathbf{A}(\mathcal{R}_0(p),\Delta)>\mathbf{B}(\mathcal{R}_0(p),\Delta).
\label{eq:AB-inequality}
\end{equation}

\begin{remark}
For an arbitrary consecutive CA pair $M$ and $N$, $G(N)$ can be greater
or less than $G(M)$, which lets $\mathbf{A}(\mathcal{R}_0(p),\Delta)$ and $\mathbf{B}(\mathcal{R}_0(p),\Delta)$
fall on either side of one another. The direction of the inequality in
Lemma~\ref{lem:master-tracking} is thus a
consequence of $N$ crossing the Robin threshold, not a generic feature of
consecutive CA numbers.
\end{remark}

\begin{lemma}[Direction of $\mathbf{A}$ and $\mathbf{B}$]\label{lem:direction}
\emergencystretch=3em For fixed $M>5040$, $\mathcal{R}_0(p)$ and $v$, let $\mathbf{A}(\mathcal{R}_0(p),y)$ and $\mathbf{B}(\mathcal{R}_0(p),y)$ be given by \eqref{eq:AB-definition} with the real variable $y$ in place of $\Delta$. Then $\mathbf{A}(\mathcal{R}_0(p),y)$ is strictly decreasing in $y$ for all real $y$, and $\mathbf{B}(\mathcal{R}_0(p),y)$ is strictly increasing in $y$ on its domain
\[
y>-(\ln M+\mathcal{R}_0(p)).
\]
\end{lemma}

\begin{proof}
By \eqref{eq:AB-definition}, $\mathbf{A}(\mathcal{R}_0(p),y)=1+1\big/\sum_{k=1}^{v+1}e^{k(\mathcal{R}_0(p)+y)}$. Each term $e^{k(\mathcal{R}_0(p)+y)}$ is positive and strictly increasing in $y$. The sum is therefore positive and strictly increasing in $y$, which makes $\mathbf{A}(\mathcal{R}_0(p),y)$ strictly decreasing in $y$ for all real $y$. For $\mathbf{B}$, $\ln\ln M>0$ because $M>5040$, and $\ln M+\mathcal{R}_0(p)+y$ is strictly increasing in $y$ and positive throughout the domain. Since the logarithm is strictly increasing, $\mathbf{B}(\mathcal{R}_0(p),y)$ is strictly increasing in $y$ on its domain.
\end{proof}

\subsection{The Single-Prime Case}

\begin{remark}
Lemma~\ref{lem:master-tracking} says nothing about the non-trivial zeros of $\zeta(s)$. Rather, the lemma is a constraint on the prime $p$. This constraint carries over to the von Mangoldt expansion applied to $\ln p$, because $\ln p=\mathcal{R}_0(p)+\Delta(p,\{\delta_\rho\})$ by \eqref{eq:R0-Delta-sum}, and this turns Lemma~\ref{lem:master-tracking} into the inequality \eqref{eq:AB-inequality} in $\Delta$. In this way the non-trivial zeros enter the constraint. Consider $\mathbf{A}(\mathcal{R}_0(p),y)$ and $\mathbf{B}(\mathcal{R}_0(p),y)$ for real $y$, as in Lemma~\ref{lem:direction}. $\mathbf{A}$ is positive for all $y$. $\mathbf{B}$ is negative when $0<\ln M+\mathcal{R}_0(p)+y<1$, and $\mathbf{B}$ exceeds $\mathbf{A}$ for all sufficiently large $y$, because $\mathbf{B}$ grows without bound while $\mathbf{A}$ tends to $1$. Hence, by Lemma~\ref{lem:direction}, $\mathbf{A}$ and $\mathbf{B}$ are equal at exactly one value of $y$. As we will show, in the convergent case, if that value of $y$ is at most $0$, then $\Delta(p,\{\delta_\rho\})$ must be negative. An arithmetic fact about the least CA counterexample $N$ and its CA predecessor $M$ leads to information on how the off-line zeros collectively behave when $N$ crosses the Robin threshold.
\end{remark}

\begin{theorem}[Baseline Inequality for $\Delta < 0$]
\label{thm:main_spectral_constraint}
Let $N=M\cdot p$ be the least CA counterexample to Robin's Inequality
in the single-prime case, with $M$ its CA predecessor and $v$ the $p$-adic valuation of $M$. Let $\Delta = \Delta(p,\{\delta_\rho\})$ denote the off-line distortion
in the von Mangoldt explicit formula expansion for $\ln p$. In the convergent case, assume $\ln M + \mathcal{R}_0(p) > 0$ (needed for $\mathbf{B}(\mathcal{R}_0(p),0)$ to be well-defined), and suppose
\begin{equation}
1+\frac{1}{e^{\mathcal{R}_0(p)}+e^{2\mathcal{R}_0(p)}+\cdots+e^{(v+1)\mathcal{R}_0(p)}} \;\le\; \frac{\ln(\ln M+\mathcal{R}_0(p))}{\ln\ln M},
\label{eq:A0-leq-B0}
\end{equation}
or equivalently $\mathbf{A}(\mathcal{R}_0(p),0)\le\mathbf{B}(\mathcal{R}_0(p),0)$. Then
\begin{equation}
\Delta < 0.
\end{equation}
\end{theorem}

\begin{proof}
Suppose, for contradiction, $\Delta\ge0$. Since $\ln M+\mathcal{R}_0(p)>0$, both $0$ and $\Delta$ lie in the domain of $\mathbf{B}$ in Lemma~\ref{lem:direction}. That lemma then gives $\mathbf{A}(\mathcal{R}_0(p),\Delta)\le\mathbf{A}(\mathcal{R}_0(p),0)$ and $\mathbf{B}(\mathcal{R}_0(p),0)\le\mathbf{B}(\mathcal{R}_0(p),\Delta)$. Together with the inequality \eqref{eq:A0-leq-B0},
\begin{equation}
\mathbf{A}(\mathcal{R}_0(p),\Delta) \;\le\; \mathbf{A}(\mathcal{R}_0(p),0) \;\le\; \mathbf{B}(\mathcal{R}_0(p),0) \;\le\; \mathbf{B}(\mathcal{R}_0(p),\Delta),
\end{equation}
contradicting Lemma~\ref{lem:master-tracking} in the form \eqref{eq:AB-inequality}. Hence $\Delta < 0$.
\end{proof}

\subsection{The Semiprime Case}
The single-prime analysis of Theorem~\ref{thm:main_spectral_constraint} rests on the
inequality of Lemma~\ref{lem:master-tracking} together with Lemma~\ref{lem:direction}. We now
generalize that argument to cover the semiprime case $N=M\cdot pq$, where
two off-line distortions $\Delta_p,\Delta_q$ must be tracked jointly.

\begin{definition}[Semiprime notation $\mathbf{A}_p, \mathbf{A}_q, \mathbf{B}_{pq}$]
\label{def:Ap-Aq-semiprime}
In the semiprime case $N=M\cdot pq$, with $M$ the CA predecessor of $N$ and $p<q$ distinct primes, let $v_p$ and $v_q$ be the $p$-adic and $q$-adic valuations of $M$, and let $\mathcal{R}_0(p)$ and $\mathcal{R}_0(q)$ be the baseline contributions of Section~\ref{sec:baseline} and $\Delta_p$ and $\Delta_q$ the off-line distortions of Section~\ref{sec:offline-distortion}, for the pairs $(p,v_p)$ and $(q,v_q)$. By \eqref{eq:R0-Delta-sum},
\[
\Delta_p=\ln p-\mathcal{R}_0(p), \qquad \Delta_q=\ln q-\mathcal{R}_0(q).
\]
Set
$\Delta_{pq}\coloneqq\Delta_p+\Delta_q$. Define $\mathbf{A}_p(\mathcal{R}_0(p),\Delta_p)$ and $\mathbf{A}_q(\mathcal{R}_0(q),\Delta_q)$
by the same formula as $\mathbf{A}(\mathcal{R}_0(p),\Delta)$ in \eqref{eq:AB-definition}, applied
separately to each prime's own baseline contribution and valuation
\begin{equation}
\mathbf{A}_p(\mathcal{R}_0(p),\Delta_p) \coloneqq 1+\frac{1}{\sum_{k=1}^{v_p+1}e^{k(\mathcal{R}_0(p)+\Delta_p)}},
\qquad
\mathbf{A}_q(\mathcal{R}_0(q),\Delta_q) \coloneqq 1+\frac{1}{\sum_{k=1}^{v_q+1}e^{k(\mathcal{R}_0(q)+\Delta_q)}}.
\label{eq:Ap-Aq-sum-form}
\end{equation}
Since the added primes contribute $\mathcal{R}_0(p) + \mathcal{R}_0(q) + \Delta_{pq}$ to $\ln N$, the semiprime case uses $\mathbf{B}_{pq}$
with $\mathcal{R}_0(p)$ and $\mathcal{R}_0(q)$ combined, in place of the single $\mathcal{R}_0(p)$ of
\eqref{eq:AB-definition}
\begin{equation}
\mathbf{B}_{pq}(\mathcal{R}_0(p),\mathcal{R}_0(q),\Delta_{pq}) \coloneqq \frac{\ln(\ln M+\mathcal{R}_0(p)+\mathcal{R}_0(q)+\Delta_{pq})}{\ln\ln M}.
\end{equation}
\end{definition}

\begin{lemma}[Semiprime Master Tracking Inequality]
\label{lem:semiprime-master-tracking}
In the setting of Definition~\ref{def:Ap-Aq-semiprime}, with $N$ the least CA counterexample to Robin's Inequality,
\begin{equation}
\mathbf{A}_p(\mathcal{R}_0(p),\Delta_p)\,\mathbf{A}_q(\mathcal{R}_0(q),\Delta_q) \;>\; \mathbf{B}_{pq}(\mathcal{R}_0(p),\mathcal{R}_0(q),\Delta_{pq}).
\end{equation}
\end{lemma}

\begin{proof}
By \eqref{eq:R0-Delta-sum}, applied independently to the pair $(p,v_p)$ and to
the pair $(q,v_q)$,
\begin{equation}
\ln p = \mathcal{R}_0(p)+\Delta_p, \qquad \ln q = \mathcal{R}_0(q)+\Delta_q,
\end{equation}
each identity depending only on its own prime and valuation. Summing gives
\[
\ln p+\ln q = \big(\mathcal{R}_0(p)+\mathcal{R}_0(q)\big)+\Delta_{pq}.
\]
Hence
\begin{equation}
\ln N = \ln M+\ln p+\ln q = \ln M+\mathcal{R}_0(p)+\mathcal{R}_0(q)+\Delta_{pq}.
\end{equation}

From $G(M)<e^\gamma$ and Definition~\ref{def:G}, $\sigma(M)/M<e^\gamma\ln\ln M$. Let $M'$ denote the part of $M$ coprime to both $p$ and $q$, so that
$M=p^{v_p}q^{v_q}M'$.
Then $N=Mpq=p^{v_p+1}q^{v_q+1}M'$. Since $\sigma$ is multiplicative over the coprime
factors $p^{v_p+1}$, $q^{v_q+1}$, $M'$, the factor of $M'$ cancels and
\begin{equation}
\frac{\sigma(N)/N}{\sigma(M)/M} =
\frac{\sigma(p^{v_p+1})/p^{v_p+1}}{\sigma(p^{v_p})/p^{v_p}} \cdot
\frac{\sigma(q^{v_q+1})/q^{v_q+1}}{\sigma(q^{v_q})/q^{v_q}}.
\end{equation}
By the same geometric-sum computation as in
\eqref{eq:abundancy-scaling} and the display after it, with $v_p$ and $v_q$ in place of $v$, each factor equals $1+1/(p+p^2+\cdots+p^{v_p+1})$ and $1+1/(q+q^2+\cdots+q^{v_q+1})$ respectively. Since $p=e^{\mathcal{R}_0(p)+\Delta_p}$ and $q=e^{\mathcal{R}_0(q)+\Delta_q}$, by Definition~\ref{def:Ap-Aq-semiprime} they equal
$\mathbf{A}_p(\mathcal{R}_0(p),\Delta_p)$ and $\mathbf{A}_q(\mathcal{R}_0(q),\Delta_q)$, which gives
\begin{equation}
\frac{\sigma(N)}{N} = \frac{\sigma(M)}{M}\cdot\mathbf{A}_p(\mathcal{R}_0(p),\Delta_p)\,\mathbf{A}_q(\mathcal{R}_0(q),\Delta_q).
\end{equation}

Since $\mathbf{A}_p$ and $\mathbf{A}_q$ are positive, multiplying the bound on
$\sigma(M)/M$ by their product gives
\begin{equation}
\frac{\sigma(N)}{N} < e^\gamma\ln\ln M\cdot\mathbf{A}_p(\mathcal{R}_0(p),\Delta_p)\,\mathbf{A}_q(\mathcal{R}_0(q),\Delta_q).
\end{equation}
Conversely, the defining condition $G(N)\ge e^\gamma$ requires
\begin{equation}
\frac{\sigma(N)}{N} \ge e^\gamma\ln\ln N =
e^\gamma\ln\big(\ln M+\mathcal{R}_0(p)+\mathcal{R}_0(q)+\Delta_{pq}\big).
\end{equation}
Chaining these bounds on $\sigma(N)/N$ and dividing by $e^\gamma\ln\ln M>0$ isolates
the product, giving
\begin{equation}
\frac{\ln\big(\ln M+\mathcal{R}_0(p)+\mathcal{R}_0(q)+\Delta_{pq}\big)}{\ln\ln M}
< \mathbf{A}_p(\mathcal{R}_0(p),\Delta_p)\,\mathbf{A}_q(\mathcal{R}_0(q),\Delta_q).
\end{equation}
\end{proof}

\begin{theorem}[Baseline Inequality for the Semiprime Case]
\label{thm:joint-partial-necessity}
Let $N = M\cdot pq$ be the least CA counterexample to Robin's
Inequality, with $M$ its CA predecessor and $p<q$ distinct primes.  Let $\Delta_p$, $\Delta_q$, $\Delta_{pq}$ and $\mathbf{A}_p$, $\mathbf{A}_q$, $\mathbf{B}_{pq}$ be as in Definition~\ref{def:Ap-Aq-semiprime}.
In the convergent case, assume $\ln M+\mathcal{R}_0(p)+\mathcal{R}_0(q)>0$ (needed for $\mathbf{B}_{pq}(\mathcal{R}_0(p),\mathcal{R}_0(q),0)$ to be well-defined), and suppose
\begin{equation}
\left(1+\frac{1}{\sum_{k=1}^{v_p+1}e^{k\mathcal{R}_0(p)}}\right)\left(1+\frac{1}{\sum_{k=1}^{v_q+1}e^{k\mathcal{R}_0(q)}}\right) \le \frac{\ln(\ln M+\mathcal{R}_0(p)+\mathcal{R}_0(q))}{\ln\ln M},
\label{eq:Ap0Aq0-leq-B0}
\end{equation}
or equivalently $\mathbf{A}_p(\mathcal{R}_0(p),0)\,\mathbf{A}_q(\mathcal{R}_0(q),0)\le\mathbf{B}_{pq}(\mathcal{R}_0(p),\mathcal{R}_0(q),0)$.
Then $\Delta_p\ge0$ and $\Delta_q\ge0$ cannot both hold, that is, $\Delta_p<0$ or $\Delta_q<0$.
More precisely, whenever $\Delta_{pq}\ge0$, either
\begin{equation}
\mathbf{A}_p(\mathcal{R}_0(p),\Delta_p)\,\mathbf{A}_q(\mathcal{R}_0(q),0) \;>\; \mathbf{B}_{pq}(\mathcal{R}_0(p),\mathcal{R}_0(q),0),
\end{equation}
or
\begin{equation}
\mathbf{A}_p(\mathcal{R}_0(p),0)\,\mathbf{A}_q(\mathcal{R}_0(q),\Delta_q) \;>\; \mathbf{B}_{pq}(\mathcal{R}_0(p),\mathcal{R}_0(q),0),
\end{equation}
each implying that the corresponding distortion is strictly negative.
\end{theorem}
\begin{proof}
\emergencystretch=4em
If $\Delta_{pq}<0$, then since $\Delta_{pq}=\Delta_p+\Delta_q$, at least one of
$\Delta_p,\Delta_q$ must already be negative. The conclusion holds immediately in this case.

Assume instead $\Delta_{pq}\ge0$. By the argument of Lemma~\ref{lem:direction}, applied to the formulas of Definition~\ref{def:Ap-Aq-semiprime}, $\mathbf{A}_p$ and $\mathbf{A}_q$ are strictly decreasing in their second arguments, and $\mathbf{B}_{pq}$ is strictly increasing in its third argument over the values greater than $-(\ln M+\mathcal{R}_0(p)+\mathcal{R}_0(q))$. Since $\ln M+\mathcal{R}_0(p)+\mathcal{R}_0(q)>0$, both $0$ and $\Delta_{pq}$ are among those values. If $\Delta_p\ge0$ and $\Delta_q\ge0$ both
held, then with $\mathbf{A}_p$ and $\mathbf{A}_q$ positive and decreasing and $\mathbf{B}_{pq}$ increasing, we would have
\[
\begin{aligned}
&\mathbf{A}_p(\mathcal{R}_0(p),\Delta_p)\,\mathbf{A}_q(\mathcal{R}_0(q),\Delta_q) \\
&\le \mathbf{A}_p(\mathcal{R}_0(p),0)\,\mathbf{A}_q(\mathcal{R}_0(q),0) \\
&\le \mathbf{B}_{pq}(\mathcal{R}_0(p),\mathcal{R}_0(q),0) \\
&\le \mathbf{B}_{pq}(\mathcal{R}_0(p),\mathcal{R}_0(q),\Delta_{pq}),
\end{aligned}
\]
contradicting Lemma~\ref{lem:semiprime-master-tracking}. Thus, there are two cases, either $\Delta_p<0$ or $\Delta_q<0$.

\textit{Case $\Delta_p<0$.} Since $\Delta_{pq}\ge0$, we have $\Delta_q\ge-\Delta_p>0$. Since $\mathbf{A}_q$ is strictly decreasing and $\Delta_q\ge0$ and $\mathbf{A}_p$ is positive, Lemma~\ref{lem:semiprime-master-tracking} gives
\begin{equation}
\begin{split}
\mathbf{A}_p(\mathcal{R}_0(p),\Delta_p)\,\mathbf{A}_q(\mathcal{R}_0(q),0) &\;\ge\;
\mathbf{A}_p(\mathcal{R}_0(p),\Delta_p)\,\mathbf{A}_q(\mathcal{R}_0(q),\Delta_q) \\
&\;>\; \mathbf{B}_{pq}(\mathcal{R}_0(p),\mathcal{R}_0(q),\Delta_{pq}).
\end{split}
\label{eq:prod_bound_refined}
\end{equation}
Since $\mathbf{B}_{pq}$ is strictly increasing in its third argument and $\Delta_{pq}\ge0$,
\[
\mathbf{B}_{pq}(\mathcal{R}_0(p),\mathcal{R}_0(q),\Delta_{pq})\ge\mathbf{B}_{pq}(\mathcal{R}_0(p),\mathcal{R}_0(q),0).
\]
Combining with \eqref{eq:prod_bound_refined},
\begin{equation}
\mathbf{A}_p(\mathcal{R}_0(p),\Delta_p)\,\mathbf{A}_q(\mathcal{R}_0(q),0) \;>\; \mathbf{B}_{pq}(\mathcal{R}_0(p),\mathcal{R}_0(q),\Delta_{pq}) \;\ge\; \mathbf{B}_{pq}(\mathcal{R}_0(p),\mathcal{R}_0(q),0).
\label{eq:chain-step}
\end{equation}
Were $\Delta_p\ge0$, then since $\mathbf{A}_p$ decreases and by \eqref{eq:Ap0Aq0-leq-B0} divided by $\mathbf{A}_q(\mathcal{R}_0(q),0)$, we would have
\begin{equation}
\mathbf{A}_p(\mathcal{R}_0(p),\Delta_p)\;\le\;\mathbf{A}_p(\mathcal{R}_0(p),0)\;\le\;\mathbf{B}_{pq}(\mathcal{R}_0(p),\mathcal{R}_0(q),0)/\mathbf{A}_q(\mathcal{R}_0(q),0),
\end{equation}
that is,
\[
\mathbf{A}_p(\mathcal{R}_0(p),\Delta_p)\,\mathbf{A}_q(\mathcal{R}_0(q),0)\le\mathbf{B}_{pq}(\mathcal{R}_0(p),\mathcal{R}_0(q),0),
\]
contradicting \eqref{eq:chain-step}. Hence this inequality implies $\Delta_p<0$.

\textit{Case $\Delta_q<0$.} Since $\Delta_{pq}\ge0$, we have $\Delta_p\ge-\Delta_q>0$. The identical argument with the roles of
$p$ and $q$ exchanged, using $\Delta_p\ge0$ in place of $\Delta_q\ge0$, gives the second inequality of the statement and implies $\Delta_q<0$.

In either case the conclusion holds.
\end{proof}

\begin{remark}[The residual mixed-sign region]
Theorem~\ref{thm:joint-partial-necessity} excludes the closed quadrant $\Delta_p\ge0,\Delta_q\ge0$. At any other point $(\Delta_p,\Delta_q)$, either $\Delta_{pq}<0$, or one distortion is negative and the other is positive enough to make their sum at least zero. The points of the second kind form the residual region. The following table states the status of the cases in which a distortion is zero and the cases in which the distortions have opposite signs.
\begin{center}
\renewcommand{\arraystretch}{1.25}
\begin{tabular}{@{}>{\raggedright\arraybackslash}p{0.52\textwidth}>{\raggedright\arraybackslash}p{0.42\textwidth}@{}}
\hline
Case & Status \\
\hline
$\Delta_p=\Delta_q=0$ & excluded by Theorem~\ref{thm:joint-partial-necessity} \\
one distortion $0$, the other $>0$ & excluded by Theorem~\ref{thm:joint-partial-necessity} \\
one distortion $0$, the other $<0$ & $\Delta_{pq}<0$ \\
opposite signs, the negative distortion larger in magnitude & $\Delta_{pq}<0$ \\
opposite signs, equal in magnitude & $\Delta_{pq}=0$, in the residual region \\
opposite signs, the positive distortion larger in magnitude & $\Delta_{pq}>0$, in the residual region \\
\hline
\end{tabular}
\end{center}
\end{remark}

\section{Summary}
For the single-prime case, Lemma~\ref{lem:master-tracking} gives the inequality \eqref{eq:AB-inequality} for the
least CA counterexample. In the convergent case, Theorem~\ref{thm:main_spectral_constraint} gives $\Delta<0$ whenever the inequality \eqref{eq:A0-leq-B0} holds. Explicitly, that inequality is
\begin{equation}
1+\frac{1}{\sum_{k=1}^{v+1}e^{k\mathcal{R}_0(p)}} \;\le\; \frac{\ln(\ln M+\mathcal{R}_0(p))}{\ln\ln M}.
\end{equation}
The sign combinations for the pair $(\mathcal{R}_0(p),\Delta)$ restrict where the inequality can hold. Because $\ln p$ is the sum of $\mathcal{R}_0(p)$ and $\Delta$ by \eqref{eq:R0-Delta-sum} and $\ln p>0$ since
$p\ge2$, the sign of $\mathcal{R}_0(p)$ and the sign of $\Delta$ cannot be chosen
independently. Of the nine sign combinations for the pair
$(\mathcal{R}_0(p),\Delta)$, four are arithmetically impossible, since $(0,0)$,
$(-,-)$, $(-,0)$, and $(0,-)$ all imply $\mathcal{R}_0(p)+\Delta\le0$. Of the
remaining five, four have $\Delta\ge0$ (namely $(+,+)$, $(+,0)$, $(0,+)$, and $(-,+)$),
and each of these already makes the inequality \eqref{eq:A0-leq-B0} fail directly, by Theorem~\ref{thm:main_spectral_constraint}, provided $\ln M+\mathcal{R}_0(p)>0$, which is needed for $\mathbf{B}(\mathcal{R}_0(p),0)$ to be defined. This leaves one sign
combination for which the inequality \eqref{eq:A0-leq-B0} can hold at all. It is
$\mathcal{R}_0(p)>0$ and $\Delta<0$.

For the semiprime case, Theorem~\ref{thm:joint-partial-necessity}
establishes that $\Delta_p<0$ or $\Delta_q<0$ whenever the inequality
\eqref{eq:Ap0Aq0-leq-B0} is satisfied. It leaves the sign of $\Delta_{pq}$ open, because the residual region, where one distortion is negative and the other is positive enough to make their sum at least zero, remains possible.

In summary, this paper does not determine whether a least CA
counterexample exists. It identifies, in the single-prime case, the
inequality \eqref{eq:A0-leq-B0} for the least CA counterexample
and shows that the off-line distortion is negative if that inequality holds. In the semiprime case, the analogous inequality gives only that at least one of the off-line distortions of the von Mangoldt expansions for $\ln p$ and $\ln q$ is negative. The main results of this paper treat only the convergent case. Proving the Alaoglu--Erdős conjecture
would preempt the semiprime case entirely, since it would imply that the quotient
of consecutive CA numbers is always prime.

\end{document}